\documentclass[a4paper,11pt]{amsart}

\usepackage{amssymb,verbatim,hyperref}

\theoremstyle{plain}
\newtheorem{lem}{Lemma}[section]
\newtheorem{prop}[lem]{Proposition}
\newtheorem{cor}[lem]{Corollary}
\newtheorem{thm}[lem]{Theorem}

\theoremstyle{definition}

\theoremstyle{remark}
\newtheorem*{note}{Note}
\newtheorem*{eg*}{Example}

\newtheorem*{rem}{Remark}

\newcommand{\func}{\operatorname}
\newcommand{\Ext}{\func{Ext}}
\newcommand{\soc}{\func{soc}}
\newcommand{\wt}{\widetilde}
\newcommand{\wb}{\overline}
\newcommand{\wh}{\widehat}
\newcommand{\sym}{\mathfrak{S}}
\newcommand{\up}[1]{\text{$\uparrow^{#1}$}}
\newcommand{\down}[1]{\text{$\downarrow_{#1}$}}

\newcommand{\pty}[1]{\sigma_p({#1})}
\newcommand{\sgn}{\func{sgn}}
\newcommand{\rad}{\func{rad}}
\newcommand{\End}{\func{End}}
\newcommand{\Hom}{\func{Hom}}
\newcommand{\bd}{\bullet}

\begin{document}

\title[$v$-decomposition numbers and symmetric groups]{Parities of $v$-decomposition numbers and \\an application to symmetric group algebras}
\author{Kai Meng Tan}
\address{Department of Mathematics, National University of
Singapore, 2, Science Drive 2, Singapore 117543.}
\email{tankm@nus.edu.sg}
\date{September 2006}
\thanks{The author, supported by Academic Research Fund R-146-000-043-112 of
National University of Singapore, thanks Bernard Leclerc for his suggestion of the proof of Theorem \ref{T:appendix}, and Hung Yean Loke for directing him to \cite{IM} in which the fomula for the length of a general element of an affine Weyl group can be found.}

\begin{abstract}
We prove that the $v$-decomposition number $d_{\lambda\mu}(v)$ is an even or odd polynomial according to whether the partitions $\lambda$ and $\mu$ have the same relative sign (or parity) or not.  We then use this result to verify Martin's conjecture for weight $3$ blocks of symmetric group algebras --- that these blocks have the property that their projective (indecomposable) modules have a common radical length $7$.
\end{abstract}

\maketitle

\section{Introduction}

Throughout this paper, let $v$ be an indeterminate, and let $e$ be an integer greater than
1.  The Fock space representation $\mathcal{F}$ of
$U_v(\widehat{\mathfrak{sl}}_e)$, as a $\mathbb{C}(v)$-vector space,
has two distinguished bases, the standard basis $\{ s(\lambda) \mid
\lambda \in \mathcal{P} \}$ and the canonical basis $\{ G(\lambda)
\mid \lambda \in \mathcal{P} \}$, both being indexed by the set
$\mathcal{P}$ of all partitions of non-negative integers.  The
$v$-decomposition number $d_{\lambda\mu}(v) \in \mathbb{C}(v)$ is
the coefficient of $s(\lambda)$ when the canonical basis element
$G(\mu)$ is expressed in terms of the standard basis elements, i.e.\
$$
G(\mu) = \sum_{\lambda \in \mathcal{P}} d_{\lambda\mu}(v)
s(\lambda).
$$

Varagnolo and Vasserot \cite{VV} showed that the $v$-decomposition numbers are parabolic Kazhdan-Lusztig's polynomials.  As an immediate consequence, a non-zero $v$-decomposition number is a sum of monomials in $v$, either all of which are of odd degree or of even degree.  In the first part of this paper, we provide a combinatorial criteria, in terms of relative signs (or parities) of the partitions $\lambda$ and $\mu$, that determines exactly which of the two possibilities occurs.

In the second part of this paper, we investigate a conjecture of Martin concerning the blocks of symmetric group algebras with Abelian defect groups.  A block of a symmetric group algebra in positive characteristic $p$ is parametrised by a pair $(\kappa, w)$, where $\kappa$ is a $p$-core partition and $w$ is a non-negative integer called its ($p$-)weight of the block.  The defect group of such a block is Abelian if and only if $w < p$.  Martin's conjecture \cite{Martin} asserts that in such a case, the projective (indecomposable) modules of this block have a common radical length $2w+1$.  This conjecture have been verified for $w \leq 2$, and for the case of $w=3$, Martin and the author \cite{mt1,mt2} showed that certain weight $3$ blocks indeed have such a property, and obtained some sufficient conditions for the inheritance of such a property for one weight $3$ block from another with which they form a $[3:1]$- or $[3:2]$-pair.  One of these sufficient conditions is the assumption that the pair of blocks concerned have bipartite Ext-quivers.  This condition is recently shown to hold for all weight $3$ blocks in \cite{ft1}, with the bipartition being given by the relative signs (or parities) of partitions labelling the simple modules.  Our result on the parities of the $v$-decomposition numbers can therefore be applied to obtain information on the module structures of some important modules.  This new ingredient proves sufficient to provide a complete verification of Martin's conjecture for $w=3$.

We now indicate the layout of this paper.  In the remainder of this section, we introduce some non-standard notation and give a summary of the combinatorics of partitions which we shall require.  In Section \ref{S:Fock}, we review the theory of $v$-decomposition numbers, state the first main result of this paper and prove it while assuming Proposition \ref{P:mainprop}.  Section \ref{S:W-action} is devoted entirely to the proof of Proposition \ref{P:mainprop}.  The remainder of the paper deals with the representation theory of the symmetric groups.  We give a brief account of the general theory in Section \ref{S:symgp}, while in Section \ref{S:wt3}, we specialise to weight $3$ blocks.  We also state the second main result of this paper in Section \ref{S:wt3} and prove it while assuming Proposition \ref{P:32}.  Section \ref{S:32} is devoted entirely to the proof of Proposition \ref{P:32}.

\subsection{Notation}

Given a ring $R$, a simple left $R$-module $S$ and an arbitrary left $R$-module $M$ with a composition series, we write $[M:S]$ for the multiplicity of $S$ as a composition factor of $M$.

\subsection{Partitions}

A partition $\lambda
=(\lambda_{1},\lambda_{2},\dotsc)$ is a weakly decreasing sequence of non-negative integers, where for sufficiently large $i$, $\lambda_i =0$.  If $\sum_i \lambda_i = n$, we say $\lambda$ is a partition of $n$.  The length of $\lambda$, denoted $l(\lambda)$, equals $\max (i \mid \lambda_i >0)$.
Denote the set of partitions of
${n}$ by $\mathcal{P}_{n}$, and let
$\mathcal{P}=\bigcup_{n} \mathcal{P}_{n}$ be the set of all partitions.

A strictly decreasing sequence $\beta= (\beta_1,\beta_2,\dotsc, \beta_s)$ of non-negative integers is a sequence of $\beta$-numbers for $\lambda$ if $s \geq l(\lambda)$, and $\beta_i = \lambda_i + s - i$ for all $1 \leq i \leq s$.  Every strictly decreasing sequence of non-negative integers is a sequence of $\beta$-numbers for a unique partition.

The James $e$-abacus has $e$ vertical runners, labelled $0,1,\dotsc, e-1$.  Its positions are labelled from left to right, and top down, starting from $0$.  The partition $\lambda$ may be displayed on the abacus as follows:  if $\beta =(\beta_1,\beta_2,\dotsc, \beta_s)$ is a sequence of $\beta$-numbers for $\lambda$, then we place a bead at position $\beta_i$ for each $i$.  This is the ($e$-)abacus display of $\lambda$ with $s$ beads.

In an abacus display of $\lambda$, moving a bead from position $a$ to a vacant position $b$, with $a>b$, corresponds to removing a (rim) hook of length $a-b$ from $\lambda$.  The number of beads crossed in so doing (i.e.\ the number of occupied positions between $b$ and $a$) is the leg-length of the hook.  The $e$-core of $\lambda$ is thus obtained when we slide the beads as far up their respective runners as possible.  The $e$-weight of $\lambda$ is the total number of times we slide the beads one position up their respective runners to obtain its $e$-core.  The relative ($e$-)sign of $\lambda$, denoted as $\sigma_e(\lambda)$, can be defined as $(-1)^t$, where $t$ is the total number of beads crossed to obtain the $e$-core (see \cite[\S2]{MO}).

The conjugate partition of $\lambda$, denoted $\lambda' = (\lambda'_1,\lambda'_2,\dotsc)$, is defined by $\lambda'_j = |\{ i \mid \lambda_i \geq j \}|$ for all $j \in \mathbb{Z}^+$.  Given an abacus display of $\lambda$, we can obtain the abacus display of $\lambda'$ by rotating the abacus of $\lambda$ through an angle of $\pi$, and read the vacant positions as occupied and the occupied positions as vacant.  Thus, $\lambda'$ has the same $e$-weight as $\lambda$, and its $e$-core is the conjugate partition of the $e$-core of $\lambda$.

The partition $\lambda$ is $e$-regular if there does not exist $i$ such that $\lambda_i = \lambda_{i+1} = \dotsb = \lambda_{i+e-1} > 0$, and is $e$-restricted if $\lambda'$ is $e$-regular.  In \cite{Mullineux}, Mullineux formulated an involution $\lambda \mapsto m(\lambda)$ on the set of $e$-regular partitions of $\mathcal{P}_n$.  This bijection plays an important role in the representation theory of Iwahori-Hecke algebras of the symmetric groups and the Fock space representation of $U_v(\widehat{\mathfrak{sl}}_e)$, which we shall describe later.  We refer the reader to \cite{Mullineux} for a combinatorial description of this bijection.

\section{The Fock space representation of $U_v(\widehat{\mathfrak{sl}}_e)$} \label{S:Fock}

In this section, we define the $v$-decomposition numbers arising from the Fock space representation of $U_v(\widehat{\mathfrak{sl}}_e)$, and briefly discuss some of the remarkable properties they enjoy.

The quantum affine algebra $U_v(\widehat{\mathfrak{sl}}_e)$ is an associative algebra over $\mathbb{C}(v)$ generated by $e_r,f_r, k_r, k_r^{-1}$ ($0 \leq r \leq e-1$), $d, d^{-1}$ subject to certain relations which we do not need here.  The Fock space representation $\mathcal{F}$ of $U_v(\widehat{\mathfrak{sl}}_e)$ is a $\mathbb{C}(v)$-vector space with basis $\{ s(\lambda) \mid \lambda \in \mathcal{P} \}$.  For our purposes, an explicit description of the actions of $e_r$ and $f_r$ will suffice.

Display a partition $\lambda$ on the $e$-abacus with $t$ beads, where $t \geq l(\lambda)$ and $e \nmid (r+t)$.  Let $i$ be the residue class of $(r+t)$ modulo $e$.  Suppose there is a bead on runner $i-1$ whose succeeding position on runner $i$ is vacant; let $\mu$ be the partition obtained when this bead is moved to its succeeding position.
Let $N_>(\lambda,\mu)$ (resp.\ $N_<(\lambda,\mu)$) be the number of beads on runner $i-1$ below (resp.\ above) the bead moved to obtained $\mu$ minus the number of beads on runner $i$ below (resp.\ above) the vacant position that becomes occupied in obtaining $\mu$.
We have
\begin{align*}
f_r( s(\lambda)) &= \sum_{\mu} v^{N_>(\lambda,\mu)} s(\mu); \\
e_r( s(\mu)) &= \sum_{\lambda} v^{-N_<(\lambda,\mu)} s(\lambda),
\end{align*}
where $\mu$ in the first sum runs over all partitions that can be obtained from $\lambda$ by moving a bead on runner $i-1$ to its vacant succeeding position on runner $i$, while $\lambda$ in the second sum runs over all partitions that can be obtained from $\mu$ by moving a bead on runner $i$ to its vacant preceding position on runner $i-1$.

Let $\left< - , - \right>$ be the inner product on $\mathcal{F}$ with respect to which $\{ s(\lambda) \mid \lambda \in \mathcal{P} \}$ is orthonormal.  With respect to this inner product, the operators $e_r$ and $f_r$ are adjoints to each other, i.e.\ $\langle f_r x, y \rangle = \langle  x, e_r y \rangle$ for all $x,y \in \mathcal{F}$.

The Fock space $\mathcal{F}$ possesses another distinguished basis $\{ G(\lambda) \mid \lambda \in \mathcal{P} \}$, called the canonical basis.  For $\lambda,\mu \in \mathcal{P}$, define the $v$-decomposition number $d_{\lambda\mu}(v) \in \mathbb{C}(v)$ to be the coefficient of $s(\lambda)$ in the expansion of $G(\mu)$, i.e.
$$
d_{\lambda\mu}(v) = \langle G(\mu), s(\lambda) \rangle.
$$

Varagnolo and Vasserot \cite{VV} showed that the $v$-decomposition numbers are parabolic Kazhdan-Lusztig's polynomials.  We give a brief account of this.

The extended affine Weyl group $W = \sym_n \ltimes \mathbb{Z}^n$ acts on $\mathbb{Z}^n$ via
$$
  \sigma (t_1,\dotsc,t_n) \cdot (a_1,\dotsc,a_n) = (t_{\sigma^{-1}(1)} + a_{\sigma^{-1}(1)},\dotsc, t_{\sigma^{-1}(n)}+ a_{\sigma^{-1}(n)}).
$$
(Here, and hereafter, $\sym_n$ denotes the symmetric group on $n$ letters.)
The set $\mathcal{A} = \{ (x_1,\dotsc, x_n) \in \mathbb{Z}^n \mid -e < x_1 \leq \dotsb \leq x_n \leq 0 \}$ is a fundamental domain of this action, and $\{ \varepsilon_i - \varepsilon_j \mid i \ne j\}$ is a root system for $W$, where $\varepsilon_i$ is the $i$-th standard basis element of $\mathbb{Z}^n$.  We take the positive roots to be $\{ \varepsilon_i - \varepsilon_j \mid i < j\}$.  For $a \in \mathbb{Z}^n$, write $w_a$ for the (unique) element in $\sym_{n} \ltimes \mathbb{Z}^n$ having the minimal length (with respect to this positive root system) such that $w_a^{-1} \cdot a \in \mathcal{A}$.

Given a partition $\lambda$ with $l(\lambda) \leq n$, write $\wh{\lambda}$ for the strictly increasing sequence of non-negative integers $(a_1,\dotsc,a_n)$ such that $(a_n,\dotsc,a_1)$ is a sequence of $\beta$-numbers of $\lambda$. Let $\mu$ be another partition with $l(\mu) \leq n$, and assume $\lambda$ and $\mu$ have the same $e$-weight and $e$-core.  Then $\wh{\lambda}$ and $\wh{\mu}$ lie in the same $W$-orbit, which intersects $\mathcal{A}$ at $\alpha$, say.  Let $w^\alpha$ be the longest element in the stabilizer of $\alpha$ (under the action of $W$).  We have

\begin{thm}[{\cite{VV}; see also \cite[Theorem 13]{Leclerc}}]
$$
d_{\lambda\mu}(v) = \sum_{y \in \sym_r} (-v)^{\ell(y)} P_{yw_{\wh{\lambda}}w^\alpha,w_{\wh{\mu}}w^\alpha}(v),
$$
where $\ell(y)$ is the length of $y$ as an element of $W$ (with respect to the chosen positive root system), and $P_{x,w}$ is the coefficient of $T_x$ in the expansion of the Kazhdan-Lusztig base element $C_w'$ of the Hecke algebra associated with $W$ (following Soergel's convention on the normalisation of the generators $T_i$ \cite{Soergel}).
\end{thm}

The upshot of this is:
\begin{cor} \label{C:par}
Keeping the above notations, we have
$$
d_{\lambda\mu}(v) \in
\begin{cases}
\mathbb{N}_0[v^2], &\text{if } (-1)^{\ell(w_{\wh{\lambda}})} = (-1)^{\ell(w_{\wh{\mu}})}; \\
v\mathbb{N}_0[v^2], &\text{if } (-1)^{\ell(w_{\wh{\lambda}})} \ne (-1)^{\ell(w_{\wh{\mu}})}.
\end{cases}
$$
\end{cor}

\begin{proof}
This follows from the fact that
$$
P_{x,w}(v) \in
\begin{cases}
\mathbb{N}_0[v^2], &\text{if } (-1)^{\ell(w)} = (-1)^{\ell(x)}; \\
v\mathbb{N}_0[v^2], &\text{if } (-1)^{\ell(w)} \ne (-1)^{\ell(x)}.
\end{cases}
$$
\end{proof}

\begin{prop} \label{P:mainprop}
Let $\lambda$ be a partition with $e$-weight $W$ and $e$-core $\kappa$, and assume that $l(\lambda) \leq n$.  Then
$$
(-1)^{\ell(w_{\wh{\lambda}})} = (-1)^{W(n-1) + \ell(w_{\wh{\kappa}})}\sigma_e(\lambda).
$$
\end{prop}

We devote the next section to the proof of Proposition \ref{P:mainprop}.  The following theorem immediately follows from Corollary \ref{C:par} and Proposition \ref{P:mainprop}:

\begin{thm} \label{T:appendix}
If $d_{\lambda\mu}(v) \ne 0$, then
$$
d_{\lambda\mu}(v) \in
\begin{cases}
\mathbb{N}_0[v^2], &\text{if } \sigma_e(\lambda) = \sigma_e(\mu); \\
v\mathbb{N}_0[v^2], &\text{if } \sigma_e(\lambda) \ne \sigma_e(\mu).
\end{cases}
$$
\end{thm}

It is proved in \cite[Proposition 2.19]{ft1} that when $p$ is an odd prime and $\lambda$ is a $p$-regular partition with $p$-weight $3$, $\pty{m(\lambda)} \ne \pty{\lambda}$.  Furthermore, it is remarked that $\pty{m(\lambda)} = (-1)^{w}\pty{\lambda}$ holds in general when $\lambda$ has $p$-weight $w$ (and is $p$-regular), by using $v$-decomposition numbers, without giving details.  Here we use Theorem \ref{T:appendix} to provide a formal proof of a more general version of this statement.

\begin{prop}
Suppose $\lambda$ is a $e$-regular partition having $e$-weight $w$.  Then $\sigma_e(m(\lambda)) = (-1)^{ew}\sigma_e(\lambda)$.
\end{prop}

\begin{proof}
Since $d_{m(\lambda)'\lambda}(v) = v^w$ \cite[Corollary 7.7]{LLT}, we have $\sigma_e(m(\lambda)')\sigma_e(\lambda) = (-1)^w$ by Theorem \ref{T:appendix}.  Note that $\sigma_e(\mu') = (-1)^{(e-1)w}\sigma_e(\mu)$ for a $e$-regular partition $\mu$ with $e$-weight $w$.  Thus, $\sigma_e(m(\lambda)) = (-1)^{ew} \sigma_e(\lambda)$.
\end{proof}

\section{The extended affine Weyl group $W = \sym_n \ltimes \mathbb{Z}^n$ action on $\mathbb{Z}^n$} \label{S:W-action}

We provide a proof of Proposition \ref{P:mainprop} in this section.  Our goal is to relate the (parity of the) length of $w_{\wh{\lambda}}$ to that of $w_{\wh{\kappa}}$, where $\kappa$ is the $e$-core of $\lambda$.

Recall the set of positive roots of $W$, and the fundamental domain $\mathcal{A}$ of the action of $W$ on $\mathbb{Z}^n$, as described in the previous section.

For $T = \sum_{i=1}^n T_i \varepsilon_i \in \mathbb{Z}^n$, and $\sigma \in \sym_n$, the length of $T\sigma$ as an element of $W$ may be calculated based on the following formula, which is a specialization of that for general extended affine Weyl groups found in Proposition 1.23 of \cite{IM}:

\begin{thm} \label{T:lengththm}
$$
\ell(T\sigma) = \sum_{\substack{i<j \\ \sigma^{-1}(i) < \sigma^{-1}(j)}} |T_i - T_j| + \sum_{\substack{i<j \\ \sigma^{-1}(i) > \sigma^{-1}(j)}} |T_i - T_j -1|.
$$
\end{thm}
Observe that the formula agrees with that for the ordinary Weyl group $\sym_n$ upon restriction.

Fix an element $a = \sum_{i=1} a_i \varepsilon_i \in \mathbb{Z}^n$, and for each $i=1,\dotsc,n$, let $t_i = \lceil a_i/e \rceil$ so that $-e < a_i - et_i \leq 0$.  Write $t = \sum_{i=1} t_i \varepsilon_i$  and $c = \sum_{i=1} c_i \varepsilon_i = -t \cdot a$.  Note that $w^{-1} \cdot a \in \mathcal{A}$ for $w \in \sym_n \ltimes \mathbb{Z}^n$ if and only if $w = t\sigma$ for some $\sigma \in \sym_n$ with $\sigma^{-1} \cdot c \in \mathcal{A}$.

A description of $w_a$ is particularly easy when $a_1 \leq a_2 \leq \dotsb \leq a_n$.

\begin{prop} \label{P:lengthcor}
If $T = \sum_{i=1}^n T_i \varepsilon_i \in \mathbb{Z}^n$ with $T_1 \leq T_2 \leq \dotsb \leq T_n$, and $\sigma \in \sym_n$, then
\begin{align*}
\ell(T \sigma) &= \sum_{i<j} (T_j - T_i) + \sum_{\substack{i<j \\ \sigma^{-1}(i) > \sigma^{-1}(j)}} 1 \\
&= \sum_{i<j} (T_j - T_i) + \ell(\sigma).
\end{align*}
In particular, if $a$, $t$ and $c$ are as defined above, with $a_1 \leq a_2 \leq \dotsb \leq a_n$, then $w_a = t\sigma_c$, where $\sigma_c$ is defined as the unique element of $\sym_n$ satisfying $\sigma_c^{-1}(i) > \sigma_c^{-1}(j)$ if and only if $c_i > c_j$, and
\begin{align*}
\ell(w_a) &= \sum_{i<j} (t_j - t_i) + l(\sigma_c) \\
&= \sum_{i<j} (t_j - t_i) + |\{ (i,j) \mid 1\leq i<j\leq n,\ c_i > c_j \}|.
\end{align*}
\end{prop}

\begin{proof}
The formula for $\ell(T \sigma)$ follows immediately from Theorem \ref{T:lengththm}.  If $a_1 \leq a_2 \leq \dotsb \leq a_n$, then necessarily $t_1 \leq t_2 \leq \dotsb \leq t_n$, so that $w_a = t\sigma_c$ where $\sigma_c$ is the shortest element of $\sym_n$ satisfying $\sigma^{-1} \cdot c \in \mathcal{A}$.  The description of $\sigma_c$ and $\ell(w_a)$ then follows.
\end{proof}

Now, assume $a_1 < a_2 < \dotsb < a_n$, and suppose $a_{s-1} < a_r -e < a_s$ for some $1\leq s \leq r \leq n$.  Let $u$ be the least index such that $t_u = t_r$.  Then $s \leq u \leq r$.  Furthermore,

\begin{lem} \hfill \label{L:tc}
\begin{enumerate}
\item $t_u = t_{u+1} = \dotsb = t_r$, and if $s < u$, then $t_s = t_{s+1} =\dotsb = t_{u-1} = t_r-1$.
\item $c_u < c_{u+1} < \dotsb < c_r$, and if $s < u$, then $c_r < c_s < c_{s+1} < \dotsb < c_{u-1}$.
\end{enumerate}
\end{lem}

\begin{proof}
Since $a_r-e < a_i < a_r$ for all $s \leq i \leq r-1$, we have $t_r - 1 \leq t_i \leq t_r$.  This yields part (1).
Part (2) then follows since $c_i = a_i - et_i$ for all $i$ (note that if $s< u$, then $c_r = a_r - et_r < a_s +e - e(t_s+1) = c_s$).
\end{proof}

Write $[s,r]$ for $(s,s+1,\dotsc, r) \in \sym_n$, and let
\begin{align*}
a' &= (a'_1,\dotsc, a'_n) = [s,r](-e\varepsilon_r) \cdot a \\
&= (a_1,\dotsc, a_{s-1}, a_r-e, a_s,\dotsc, a_{r-1}, a_{r+1},\dotsc, a_n)
\end{align*}
For each $i$, let $t'_i = \lceil a'_i/e \rceil$ and $c'_i = a'_i - et'_i$, and write $t' = (t'_1,\dotsc, t'_n)$ and  $c' = (c'_1,\dotsc, c'_n)$.  Then $c' = [s,r] \cdot c$, and

\begin{lem} \label{L:t}
$t'= t-\varepsilon_u$.
\end{lem}

\begin{proof}
This follows from part (1) of Lemma \ref{L:tc}.
\end{proof}

We wish to compare $\ell(w_a)$ and $\ell(w_{a'})$.  First, we compare $\ell(\sigma_c)$ and $\ell(\sigma_{c'})$.

\begin{lem} \label{L:c}
$\ell(\sigma_c) - \ell(\sigma_{c'}) = 2u-(r+s)$.
\end{lem}

\begin{proof}
We define a partial correspondence between the sets $X= \{ (i,j) \mid i < j,\ c_i > c_j\}$ and $Y=\{(i,j) \mid i < j,\ c'_i > c'_j \}$ as follows:
\begin{alignat*}{2}
(i,j) &\longleftrightarrow (i^+, j^+) &\qquad & (j \ne r) \\
(i,r) &\longleftrightarrow (i,s) &\qquad & (i<s),
\end{alignat*}
where
$$
i^+ =
\begin{cases}
i+1, &\text{if } s \leq i < r; \\
i, &\text{otherwise.}
\end{cases}
$$
By part (2) of Lemma \ref{L:tc}, this is actually a one-to-one correspondence between the set $X \setminus X'$ and $Y \setminus Y'$, where
\begin{align*}
X' &= \{ (i,r) \mid s \leq i < u \},\\
Y'&= \{ (s,j) \mid u < j \leq r \}.
\end{align*}
Thus, $\ell(\sigma_c) - \ell(\sigma_c') = |X|-|Y| = |X'| - |Y'| = (u-s)-(r-u) = 2u-(r+s)$.
\end{proof}

\begin{cor} \label{C:maincor}
$\ell(w_a) = \ell(w_{a'}) + 4u - n -1-r-s$.
\end{cor}

\begin{proof}
By Proposition \ref{P:lengthcor}, Lemmas \ref{L:t} and \ref{L:c}, we have
\begin{align*}
\ell(w_a) &= \sum_{i< j} (t_j - t_i) + \ell(\sigma_c) \\
&= \sum_{i<j} (t'_j - t'_i) + 2u-n-1 + \ell(\sigma_{c'}) +2u-r-s \\
&= \ell(w_{a'}) + 4u-n-1-r-s.
\end{align*}
\end{proof}

We are now ready to prove Proposition \ref{P:mainprop}.

\begin{proof}[Proof of Proposition \ref{P:mainprop}]
We prove by induction on the $e$-weight $W$ of $\lambda$, with $W=0$ being trivial.  Let $\lambda$ have positive $e$-weight $W$, and let $\mu$ be the partition obtained when one particular $e$-hook is removed from $\lambda$.  Let $\wh{\lambda} = (a_1,\dotsc, a_n)$.  Then there exist some integers $r,s$ with $1 \leq s \leq r \leq n$ such that $\wh{\mu} = [s,r](-e\varepsilon_r)\cdot \wh{\lambda}$.  Thus,
\begin{align*}
(-1)^{\ell(w_{\wh{\lambda}})} &= (-1)^{\ell(w_{\wh{\mu}})} (-1)^{(n-1) +(r-s)} \\
&= (-1)^{\ell(w_{\wh{\kappa}}) + (W-1)(n-1)}\sigma_e(\mu) (-1)^{(n-1) +(r-s)} \\
&= (-1)^{\ell(w_{\wh{\kappa}}) + W(n-1)} \sigma_e(\lambda)
\end{align*}
by Corollary \ref{C:maincor}, and induction hypothesis (note that the $e$-hook removed from $\lambda$ to obtain $\mu$ has leg-length $r-s$).
\end{proof}

\section{The symmetric group algebras} \label{S:symgp}

For the remainder of this paper, we look at the weight 3 blocks of symmetric group algebras with Abelian defect group.  Our goal is to prove Martin's conjecture that the projective (indecomposable) modules in these blocks have a common radical length $7$.  We shall see that Theorem \ref{T:appendix} will come in useful in our proof.

In this section, we give an account of the background of the representation theory of the symmetric group where the underlying field $\mathbb{F}$ has positive characteristic $p$.  Throughout, $\sym_n$ denotes the symmetric group on $n$ letters.

For each partition $\lambda = (\lambda_1,\lambda_2, \dotsc) \in \mathcal{P}_n$, one defines the Specht module $S^{\lambda}$ of $\mathbb{F}\sym_n$.
When $\lambda$ is $p$-regular, the Specht module $S^{\lambda}$ has a simple, self-dual head $D^{\lambda}$, and the modules $D^{\lambda}$, as $\lambda$ runs over all the $p$-regular partitions of $n$, give a complete list of simple modules of $\mathbb{F}\sym_n$.

The projective cover $P(D^{\mu})$ of $D^{\mu}$, where $\mu \in \mathcal{P}_n$ and is $p$-regular, has a (distinguished) filtration in which all factors are Specht modules.  The multiplicity of $S^\lambda$ as a factor in this filtration equals the composition multiplicity $[S^\lambda:D^{\mu}]$ of $D^{\mu}$ in $S^\lambda$.

Two Specht modules $S^{\lambda}$ and $S^{\mu}$ of $\mathbb{F}\sym_n$ lie in the same block if and only if $\lambda$ and $\mu$ have the same $p$-core.  As such each block $B$ of a symmetric group algebra is parametrised by an $p$-core partition $\kappa$ and a non-negative integer $w$.  A Specht module $S^{\rho}$ lies in $B$ if and only if $\rho$ has $p$-core $\kappa$ and $p$-weight $w$.  We call $\kappa$ the $p$-core of $B$, and $w$ the ($p$-)weight of $B$.  We also say $\rho$ is a partition in $B$ if $\rho$ has $p$-core $\kappa$ and $p$-weight $w$.  The defect group of $B$ is isomorphic to a Sylow $p$-subgroup of $\sym_w$.  As such, $B$ has Abelian defect group if and only if $w < p$.

\subsection{Restriction and induction} \label{SS:resind}

Let $B$ be a weight $w$
block of $\mathbb{F}\sym_n$, with $p$-core $\kappa_B$ (so $\kappa_B \in \mathcal{P}_{n-wp}$).  Fix an abacus display of $\kappa_B$, and consider the $p$-core partition $\kappa_C$ having an abacus display in which all the runners have the same number of beads as that in $\kappa_B$ except for runners $i-1$ and $i$, where respectively there are $k$ beads more and $k$ beads less than those in $\kappa_B$.  We assume that $\kappa_C \in \mathcal{P}_m$ with $m \leq n-k$, and let $C$ be the block of $\mathbb{F}\sym_{n-k}$ with $p$-core $\kappa_C$.

For every module $M$ of $B$, there is a module $M'$ of $C$ such that $M \down{C} \cong (M')^{\oplus k!}$, and for every module $N$ of $C$, there is a module $N'$ of $B$ such that $N \up{B} \cong (N')^{\oplus k!}$.  If $M$ (resp.\ $N$) is simple, then $M'$ (resp.\ $N'$) is either zero, or has a simple socle (and hence is indecomposable).  This induces, via Frobenius reciprocity, a bijection between a subset of $p$-regular partitions in $B$ with another subset of $p$-regular partitions in $C$, as follows:

\begin{thm}[See Section 11.2 of \cite{Kleshbook}]
Let $\mathcal{P}_B$ (resp. $\mathcal{P}_C$) be the set of $p$-regular partitions indexing the simple modules of $B$ (resp.\ $C$) which do not vanish upon restriction to $C$ (resp.\ induction to $B$).  Then there is a bijection $\Phi = \Phi_{B,C} : \mathcal{P}_B \to \mathcal{P}_C$ such that
$$
\soc(D^{\lambda} \down{C}) \cong (D^{\Phi(\lambda)})^{\oplus k!} \quad \text{and} \quad
\soc(D^{\Phi(\lambda)} \up{B}) \cong (D^{\lambda})^{\oplus k!}.
$$
Furthermore, $\mathcal{P}_B$, $\mathcal{P}_C$ and $\Phi$ can be combinatorially described.
\end{thm}

We have a nice description of $\End_B(D^{\Phi(\lambda)} \up{B})$ when $k=1$.

\begin{thm}[{\cite[Theorem 11.2.8(i,ii)]{Kleshbook}}] \label{T:Klesh}
Assume $k = 1$, and suppose that $\lambda \in \mathcal{P}_B$ with $[D^{\Phi(\lambda)} \up{B} : D^{\lambda}] = r$.  Then $\End_B(D^{\Phi(\lambda)} \up{B}) \cong \mathbb{F}[x]/(x^r)$, where $(x^r)$ is the ideal of the polynomial ring $\mathbb{F}[x]$ generated by $x^r$.
\end{thm}

Now suppose further that runner $i$ of the abacus display of $\kappa_B$ has exactly $k$ beads more than runner $i-1$.  Then $C$ also has weight $w$, and the abacus display of $\kappa_C$ can be obtained from that of $\kappa_B$ by interchanging the runners $i$ and $i-1$.  We say that $(B,C)$ is a $[w:k]$-pair.  In this case, $\Phi$ is a bijection between the entire sets of $p$-regular partitions in $B$ and $C$.  Furthermore, $D^{\lambda} \down{C}$ is semi-simple for many simple modules of $B$ (so $D^{\lambda} \down{C} \cong (D^{\Phi(\lambda)})^{\oplus k!}$ for these partitions) --- we call these simple modules of $B$ non-exceptional (with respect to the $[w:k]$-pair $(B,C)$), and the others exceptional.  Analogously, we define the exceptional and non-exceptional simple modules of $C$ (with respect to $(B,C)$).  It turns out that the bijection $\Phi$ sends partitions indexing exceptional (resp.\ non-exceptional) simple modules of $B$ to partitions indexing exceptional (resp.\ non-exceptional) simple modules of $C$.

A necessary and sufficient condition for the absence of exceptional simple modules of $B$ and $C$ is $w \leq k$.  In this case, $B$ and $C$ are Morita equivalent \cite{Scopes} --- we shall say that $B$ and $C$ are Scopes equivalent --- and the effect of $\Phi$ on a $p$-regular partition in $B$ is merely to interchange runners $i$ and $i-1$ of its abacus display.  It is clear that Scopes equivalence can be extended to an equivalence relation on the set of all blocks of symmetric group algebras.

\subsection{Rouquier blocks}

Let $B$ be a weight $w$ block of $\mathbb{F}\sym_n$.  We say that $B$ is Rouquier if its abacus display of its $p$-core has the following properties:  whenever runner $i$ is on the left of runner $j$, either runner $j$ has at least $w-1$ beads more than runner $i$, or runner $i$ has at least $w$ beads more than runner $j$.  It is easy to check that such a property is independent of the choice of abacus display of the $p$-core of $B$, and that the Rouquier blocks of a fixed weight form a single Scopes equivalence class.  These blocks are well understood in the Abelian defect case, by the results of \cite{ct1}.  In particular, we have the following Theorem.

\begin{thm}[{\cite[Theorem 6.4]{ct1}}] \label{T:Martin}
Suppose $B$ is a Rouquier block of weight $w$, with an Abelian defect group.  Then the projective indecomposable modules of $B$ have a common radical length $2w+1$.
\end{thm}

An arbitrary weight $w$ block can always be induced to a Rouquier block through a sequence of $[w:k]$-pairs.

\begin{lem} \label{L:sequence}
Suppose $A$ is a weight $w$ block of $\mathbb{F}\sym_n$.  Then there exists a sequence $B_0, B_1,\dotsc, B_s$ of weight $w$ blocks of symmetric group algebras such that $B_0 = A$, $B_s$ is Rouquier, and for each $1 \leq i \leq s$, $(B_{i},B_{i-1})$ is a $[w:k_i]$-pair for some $k_i \in \mathbb{Z}^+$.
\end{lem}

\begin{proof}
This is Lemma 3.1 of \cite{Fayers} in the context of the Iwahori-Hecke algebras of the symmetric groups, which is a deformation of the symmetric group algebras, and hence includes the Lemma as a special case.
\end{proof}

\subsection{Conjugate block}

Let $B$ be a block of $\mathbb{F}\sym_{n}$, with $p$-core $\kappa$.  The block of $\mathbb{F}\sym_{n}$ conjugate to $B$ is the one with $p$-core $\kappa'$.  This conjugate block, which we denote as $B'$, is Morita equivalent to $B$ via the functor $- \otimes \sgn$, where $\sgn$ is the one-dimensional sign representation of $\sym_{n}$.

If $\lambda$ is a partition in $B$, we have $S^{\lambda} \otimes \sgn \cong$ dual of $S^{\lambda'}$, and if $\lambda$ is $p$-regular, $D^{\lambda} \otimes \sgn = D^{m(\lambda)}$, where $m$ is the Mullineux map discussed earlier.  In particular, $S^{\lambda}$ has a simple socle $D^{m(\lambda')}$ when $\lambda$ is $p$-restricted.

\subsection{Connection with the Fock space}

The connection between the representation theory of the symmetric groups and the Fock space representation of $U_v(\widehat{\mathfrak{sl}}_e)$ is through the Iwahori-Hecke algebra $\mathcal{H}_n = \mathcal{H}_{\mathbb{F},q}(\sym_n)$ of the symmetric group.  As an algebra, $\mathcal{H}_n$ is generated by $T_1,T_2,\dotsc, T_{n-1}$ subject to the following relations:
\begin{alignat*}{2}
  T_iT_{i+1}T_i &= T_{i+1}T_iT_{i+1} & \quad &(1 \leq i \leq n-2);\\
  T_iT_j & = T_j T_i && (|i-j| \geq 2);\\
  (T_i -q)(T_i+1) &= 0 && (1 \leq i \leq n-1).
\end{alignat*}
Here, $q$ is an invertible element of the ground field $\mathbb{F}$, where in this subsection we allow its characteristic to be zero too.  Clearly, $\mathcal{H}_n$ is isomorphic to $\mathbb{F}\sym_n$ if $q = 1$, so that $\mathcal{H}_n$ is a deformation of $\mathbb{F}\sym_n$ in general.

Much of the representation theory of $\sym_n$ carries over to that of $\mathcal{H}_n$: for $\lambda \in \mathcal{P}_n$, we also have the Specht module $\mathcal{S}^{\lambda}$, and this has a simple head $\mathcal{D}^{\lambda}$ if $q$ is a root of unity and $\lambda$ is $e$-regular, where $e$ is the least positive integer such that $1+q+\dotsb+ q^{e-1} = 0$.  As $\lambda$ varies over the set of $e$-regular partitions in $\mathcal{P}_n$, the $\mathcal{D}^{\lambda}$'s give a complete list of non-isomorphic simple module of $\mathcal{H}_n$.  The projective cover $P(\mathcal{D}^{\lambda})$ of $\mathcal{D}^{\lambda}$ (when $\lambda$ is $e$-regular) has an analogous Specht filtation, and the blocks of $\mathcal{H}_n$ are also similarly parametrised.

When $\mathbb{F}$ has characteristic $p$, the composition multiplicity $[S^{\lambda} : D^{\mu}]$ of $D^{\mu}$ in $S^{\lambda}$ of $\mathbb{F}\sym_n$ is bounded above by the corresponding composition multiplicity $[\mathcal{S}^{\lambda} : \mathcal{D}^{\mu}]$ of $\mathcal{D}^{\mu}$ in $\mathcal{S}^{\lambda}$ of $\mathcal{H}_{\mathbb{C},q}(\sym_n)$, where $q$ is a primitive $p$-th root of unity.

Ariki \cite{Ariki} established the connection between the Fock space representation of $U_v(\widehat{\mathfrak{sl}}_e)$ and representation theory of $\mathcal{H}_n$ in characteristic zero where $q$ is a primitive $e$-th root of unity.

\begin{thm}[\cite{Ariki}]
Evaluating the $v$-decomposition number $d_{\lambda\mu}(v)$ at $v=1$ gives the composition multiplicity $[\mathcal{S}^{\lambda}:\mathcal{D}^{\mu}]$ of $\mathcal{D}^{\mu}$ in $\mathcal{S}^{\lambda}$ of $\mathcal{H}_n$ in characteristic zero.
\end{thm}

This has the following consequence:

\begin{prop} \label{P:Ariki}
Let $B$ and $C$ be blocks of $\mathcal{H}_n$ and $\mathcal{H}_{n-k}$ in characteristic zero.  Suppose that their $e$-cores, denoted $\kappa_B$ and $\kappa_C$ respectively, are of the form described in subsection \ref{SS:resind}.  Let $r$ be the residue class of $i-t$ modulo $e$, where $t$ is the number of beads used in the abacus display of $\kappa_C$ (or $\kappa_B$), and let $\lambda$ be an $e$-regular partition in $C$.  If
$$
f_r^{(k)}G(\lambda) = \sum_{\rho \in \mathcal{P}_B} a_{\rho}(v) G(\rho),$$
then
$$
P(\mathcal{D}^{\lambda}) \up{B} = \bigoplus_{\rho \in \mathcal{P}_B} P(\mathcal{D}^{\rho})^{\oplus k! a_{\rho}(1)}.
$$
\end{prop}

\begin{note}
It is sometimes possible to recover $a_{\rho}(v)$ from $a_{\rho}(1)$.  This is because $a_{\rho}(v) = a_{\rho}(v^{-1})$ and $a_{\rho}(v) \in \mathbb{N}_0[v,v^{-1}]$ (see \cite[Proposition 2.4]{Klesh}).  Thus, $a_{\rho}(v) = 0$ if and only if $a_{\rho}(1) = 0$, while $a_{\rho}(v) = 1$ if and only if $a_{\rho}(1) = 1$.  We shall use these facts later.
\end{note}

\section{Weight $3$ blocks of $\mathbb{F}\sym_n$} \label{S:wt3}

In this section, we focus our attention on weight $3$ blocks of $\mathbb{F}\sym_{n}$ with Abelian defect groups.  Thus the characteristic $p$ of $\mathbb{F}$ is assumed to be at least $5$.  These blocks enjoy many nice properties.  From now on, the $v$-decomposition numbers we are concerned with come from the Fock space representation of $U_v(\widehat{\mathfrak{sl}}_p)$ (i.e.\ $e=p$).

\begin{thm} \label{T:basicprop}
Let $B$ be a weight $3$ block of $\mathbb{F}\sym_{n}$, and let $\lambda$ and $\mu$ be partitions in $B$.
\begin{enumerate}
\item If $\mu$ is $p$-regular, then $[S^{\lambda} : D^{\mu}] = 0$ or $1$.  Furthermore, $[S^{\lambda} : D^{\mu}] = d_{\lambda\mu}(1)$.
\item If $\lambda$ and $\mu$ are $p$-regular, then $\Ext^1(D^{\lambda},D^{\mu}) = 0$ unless $\sigma_p(\lambda) \ne \sigma_p(\mu)$.  Furthermore, if $\Ext^1(D^{\lambda},D^{\mu}) \ne 0$, then $\dim_{\mathbb{F}} \Ext^1(D^{\lambda},D^{\mu}) = 1$.
\end{enumerate}
\end{thm}

\begin{proof}
Part (1) is proved by Fayers in \cite{mf3}. For part (2), the first statement is the main result of \cite{ft1}, while the second statement is proved in \cite{mr2}.
\end{proof}

\begin{cor} \label{C:monomials}
Let $\lambda$ and $\mu$ be partitions having $p$-weight $3$.  If $\mu$ is $p$-regular and $d_{\lambda\mu}(v) \ne 0$, then
$$
d_{\lambda\mu}(v) =
\begin{cases}
1, &\text{if } \lambda = \mu; \\
v, &\text{if } \lambda \notin \{ \mu, m(\mu)' \} \text{ and } \pty\lambda \ne \pty\mu; \\
v^2, &\text{if } \lambda \notin \{ \mu, m(\mu)' \} \text{ and } \pty\lambda = \pty\mu; \\
v^3, &\text{if } \lambda = m(\mu)'.
\end{cases}
$$
\end{cor}

\begin{proof}
Since $d_{\lambda\mu}(v)$ is a parabolic Kazhdan-Lusztig polynomial, it follows that $d_{\lambda\mu}(v) \in \mathbb{N}_0[v]$, so that $d_{\lambda\mu}(v)$ is a monic monomial by Theorem \ref{T:basicprop}(1).  The Corollary thus follows from Theorem 6.8 and Corollary 7.7 of \cite{LLT}, and Theorem \ref{T:appendix}.
\end{proof}

The following result is proved in \cite{ft1}.

\begin{lem}[{\cite[Proposition 2.18]{ft1}}] \label{L:pty}
Let $(B,C)$ be a $[3:k]$-pair, and let $\lambda$ be a $p$-regular partition in $B$.  Then $\sigma_p(\lambda) \ne \sigma_p(\Phi(\lambda))$ if and only if $k = 1$ and $D^{\lambda}$ is an exceptional simple module (with respect to $(B,C)$).
\end{lem}

We now state the second main result of this paper.

\begin{thm} \label{T:mainthm}
Let $B$ be a weight $3$ block of $\mathbb{F}\sym_{n}$.  Then the projective (indecomposable) modules of $B$ have a common radical length $7$.
\end{thm}

The proof of Theorem \ref{T:mainthm} relies on the following two propositions.

\begin{prop} \label{P:31}
If Theorem \ref{T:mainthm} holds for one of the blocks in a $[3:1]$-pair, then it holds for the other.
\end{prop}

\begin{proof}
Let $(B,\wt{B})$ be a $[3:1]$-pair.
From the main result of \cite{mt1} (see from page 109 onwards), we only need to show the following:
\begin{itemize}
\item whenever an exceptional simple module $D^{\alpha_j}$ of $B$ extends a non-exceptional simple module $D^{\lambda}$, then $D^{\Phi(\alpha_j)}$ does not extend $D^{\Phi(\lambda)}$;
\item whenever an exceptional simple module $D^{\Phi(\alpha_k)}$ of $\wt{B}$ extends a non-exceptional simple module $D^{\Phi(\mu)}$, then $D^{\alpha_k}$ does not extend $D^{\mu}$.
\end{itemize}
But these are immediate from Theorem \ref{T:basicprop}(2) and Lemma \ref{L:pty}.
\end{proof}

\begin{prop} \label{P:32}
Let $(B,\wt{B})$ be a $[3:2]$-pair.  If Theorem \ref{T:mainthm} holds for $B$, then it holds for $\wt{B}$.
\end{prop}

In \cite{mt2}, three sufficient conditions (Y1--Y3) for which one block in a $[3:2]$-pair may inherit Theorem \ref{T:mainthm} from the other were obtained.  However, unlike the proof of Proposition \ref{P:31}, we are unable to prove that the first of these conditions (Y1) holds in general (although we believe this to be true).  As such, we have to study $[3:2]$-pairs more carefully, which we do in the next section, to get around this and prove Proposition \ref{P:32}.

\begin{proof}[Proof of Theorem \ref{T:mainthm}]
If $B$ is Rouquier, then the Theorem holds by Theorem \ref{T:Martin}.  If $B$ is not Rouquier, then by Lemma \ref{L:sequence}, there exists a sequence of weight $3$ blocks, $B_0,B_1,\dotsc, B_s$ such that $B_0 = B$, $B_s$ is Rouquier, and for each $1 \leq i \leq s$, $(B_i, B_{i-1})$ is a $[w:k_i]$-pair for some $k_i \in \mathbb{Z}^+$.  By induction, we may assume that the Theorem holds for $B_1$.  If $k_1 \geq 3$, then $B_1$ and $B_0$ are Scopes, and hence Morita, equivalent, so that the Theorem holds for $B_0$.  If $k_1 = 1$ or $2$, then the Theorem holds for $B_0$ by Propositions \ref{P:31} and \ref{P:32} respectively.
\end{proof}

\section{$[3:2]$-pairs} \label{S:32}

In this section, $B$ is a weight 3 block of $\mathbb{F}\sym_{n}$, forming a $[3:2]$-pair with a block $\wt{B}$ of $\mathbb{F}\sym_{n-2}$.  We fix an abacus display of the $p$-core of $B$, such that by interchanging the runners $i$ and $i-1$, we obtain the abacus display of the $p$-core of $\wt{B}$.  We begin by recalling the background theory on such pairs.

Every partition in $B$, with the exception of four, has exactly two beads on runner $i$ of its abacus display which can be moved one position to the left.  The four exceptional partitions in $B$ are denoted as $\alpha$,
$\beta$, $\gamma$ and $\delta$, and the runners $i-1$ and $i$ of their respective abacus displays are as follows:
\[
  \begin{matrix}
    \begin{smallmatrix}
      i-1 & i \\
      \vdots & \vdots \\
      \bd & \bd \\
      \bd & - \\
      - & \bd \\
      - & \bd \\
      - & \bd \\
      \: & \:
    \end{smallmatrix} \\
    \alpha
  \end{matrix} \qquad
  \begin{matrix}
    \begin{smallmatrix}
      i-1 & i \\
      \vdots & \vdots \\
      \bd & \bd \\
      - & \bd \\
      \bd & - \\
      - & \bd \\
      - & \bd \\
      \: & \:
    \end{smallmatrix} \\
    \beta
   \end{matrix} \qquad
  \begin{matrix}
    \begin{smallmatrix}
      i-1 & i \\
      \vdots & \vdots \\
      \bd & \bd \\
      - & \bd \\
      - & \bd \\
      \bd & - \\
      - & \bd \\
      \: & \:
    \end{smallmatrix} \\
    \gamma
  \end{matrix} \qquad
  \begin{matrix}
    \begin{smallmatrix}
      i-1 & i \\
      \vdots & \vdots \\
      \bd & \bd \\
      - & \bd \\
      - & \bd \\
      - & \bd \\
      \bd & - \\
      \: & \:
    \end{smallmatrix} \\
    \delta
  \end{matrix}
\]

Similarly, every partition in $\wt{B}$, with the exception of four, has exactly two beads on runner $i-1$ of its abacus display which can be moved one position to the right.  The four exceptional partitions in $\wt{B}$ are denoted as $\wt{\alpha}$, $\wt{\beta}$, $\wt{\gamma}$ and $\wt{\delta}$, and the runners $i-1$ and $i$ of their respective abacus displays are as follows:
\[
  \begin{matrix}
    \begin{smallmatrix}
      i-1 & i \\
      \vdots & \vdots \\
      \bd & \bd \\
      \bd & - \\
      \bd & - \\
      \bd & - \\
      - & \bd \\
      \: & \:
    \end{smallmatrix} \\
    \wt{\alpha}
  \end{matrix} \qquad
  \begin{matrix}
    \begin{smallmatrix}
      i-1 & i \\
      \vdots & \vdots \\
      \bd & \bd \\
      \bd & - \\
      \bd & - \\
      - & \bd \\
      \bd & - \\
      \: & \:
    \end{smallmatrix} \\
    \wt{\beta}
  \end{matrix} \qquad
  \begin{matrix}
    \begin{smallmatrix}
      i-1 & i \\
      \vdots & \vdots \\
      \bd & \bd \\
      \bd & - \\
      - & \bd \\
      \bd & - \\
      \bd & - \\
      \: & \:
    \end{smallmatrix} \\
    \wt{\gamma}
  \end{matrix} \qquad
  \begin{matrix}
    \begin{smallmatrix}
      i-1 & i \\
      \vdots & \vdots \\
      \bd & \bd \\
      - & \bd \\
      \bd & - \\
      \bd & - \\
      \bd & - \\
      \: & \:
    \end{smallmatrix} \\
    \wt{\delta}
  \end{matrix}
\]

We note that the partitions \( \alpha \) and \( \wt{\alpha} \) are always \( p
\)-regular, while $\delta$ and $\wt{\delta}$ are always $p$-restricted.
In fact, $\alpha = m(\delta')$ and $\wt{\alpha} = m(\wt{\delta}')$.

The map $\Phi = \Phi_{B,\wt{B}}$ has the following effect on the exceptional partitions (when they are $p$-regular):
$$
\Phi(\alpha) = \wt{\alpha};\
\Phi(\beta) = \wt{\delta};\
\Phi(\gamma) = \wt{\gamma};\
\Phi(\delta) = \wt{\beta}.
$$
Furthermore, if $\beta$ is $p$-regular, then $[S^{\gamma} : D^{\beta}] = [S^{\delta} : D^{\beta}] = 1$; if $\gamma$ is $p$-regular, then $[S^{\delta}: D^{\gamma}] = 1 = [S^{\wt{\delta}}: D^{\wt{\gamma}}]$; if $\delta$ is $p$-regular, then $[S^{\wt{\gamma}} : D^{\wt{\beta}}] = [S^{\wt{\delta}} : D^{\wt{\beta}}] = 1$.

We note the following, which can be easily verified:
$$
\pty\alpha = \pty{\wt{\alpha}} =
\pty\gamma = \pty{\wt{\gamma}} \ne
\pty\beta = \pty{\wt{\beta}} =
\pty\delta = \pty{\wt{\delta}}.
$$

The four exceptional partitions in $\wt{B}$ have another characterisation: they are the partitions where there is a unique bead on runner $i$ of their respective abacus displays which can be moved one position to its left.  The partition so obtained is always the same, and we denote it as $\check{\alpha}$.  This partition has weight $0$ and is the unique partition in the block $\check{B}$ of $\mathbb{F}\sym_{n-3}$.  Thus $S^{\check{\alpha}} = D^{\check{\alpha}}$ is simple and projective.  By the ordinary branching rule, we have
$$
S^{\wt{\alpha}} \down{\check{B}} \cong S^{\wt{\beta}} \down{\check{B}} \cong S^{\wt{\gamma}} \down{\check{B}} \cong S^{\wt{\delta}} \down{\check{B}} \cong S^{\check{\alpha}},
$$
and $S^{\wt{\lambda}} \down{\check{B}} = 0$ for all $\wt{\lambda} \notin \{ \wt{\alpha}, \wt{\beta}, \wt{\gamma}, \wt{\delta} \}$, and hence
$D^{\wt{\alpha}} \down{\check{B}} = D^{\check{\alpha}}$ while $D^{\wt{\lambda}} \down{\check{B}} = 0$ for all $\wt{\lambda} \ne \wt{\alpha}$.  This gives us $[S^{\wt{\beta}}:D^{\wt{\alpha}}] = [S^{\wt{\gamma}}:D^{\wt{\alpha}}] = [S^{\wt{\delta}}:D^{\wt{\alpha}}] = 1$.  Furthermore, the module $D^{\check{\alpha}} \up{\check{B}}$ is projective and has a simple head $D^{\wt{\alpha}}$ by Frobenius reciprocity, and a Specht filtration filtered by $S^{\wt{\alpha}}$, $S^{\wt{\beta}}$, $S^{\wt{\gamma}}$ and $S^{\wt{\delta}}$ by the ordinary branching rule.  Thus $D^{\check{\alpha}} \up{\check{B}} \cong P(D^{\wt{\alpha}})$ and $[P(D^{\wt{\alpha}}):D^{\wt{\alpha}}] = 4$.

Entirely analogous results also hold for the exceptional partitions in $B$; i.e.\ we also have $[S^{\beta}:D^{\alpha}] = [S^{\gamma}:D^{\alpha}] = [S^{\delta}:D^{\alpha}] = 1$ and $[P(D^{\alpha}):D^{\alpha}] = 4$.

There is exactly one exceptional simple module of $\wt{B}$ with respect to $(B,\wt{B})$, namely $D^{\wt{\alpha}}$ (thus $D^{\alpha}$ is the unique exceptional simple module of $B$).

We can obtain further information on $B$ and $\wt{B}$ by considering their conjugate blocks $B'$ and $\wt{B}'$.  The latter also form a $[3:2]$-pair, and denoting the exceptional partitions in $B'$ and $\wt{B}'$ as $\alpha^c, \beta^c, \gamma^c, \delta^c$ and $\wt{\alpha}^c, \wt{\beta}^c, \wt{\gamma}^c, \wt{\delta}^c$ respectively, we have
\begin{align*}
\alpha^c = \delta',\ \beta^c = \gamma',\ \gamma^c = \beta',\ \delta^c = \alpha';\\
\wt{\alpha}^c = \wt{\delta}',\ \wt{\beta}^c = \wt{\gamma}',\ \wt{\gamma}^c = \wt{\beta}',\ \wt{\delta}^c = \wt{\alpha}'.
\end{align*}
Thus, $m(\alpha) = \alpha^c$ and $m(\wt{\alpha}) = \wt{\alpha}^c$.  As such we deduce the following:

\begin{lem} \label{L:conjugate}
We have
\begin{enumerate}
\item $\alpha$ is $p$-restricted if and only if $\wt\gamma$ is, in which case $[S^{\wt{\alpha}} : D^{m(\wt{\gamma}')}] = [S^{\wt{\beta}} : D^{m(\wt{\gamma}')}] = 1$;
\item $\beta$ is $p$-restricted if and only if $\wt\beta$ is, in which case $[S^{\alpha} : D^{m(\beta')}] = 1 = [S^{\wt{\alpha}} : D^{m(\wt{\beta}')}]$;
\item $\gamma$ is $p$-restricted if and only if $\wt\alpha$ is, in which case $[S^{\alpha} : D^{m(\gamma')}] = [S^{\beta} : D^{m(\gamma')}] = 1$;
\end{enumerate}
\end{lem}

\begin{proof}
We show part (2); parts (1) and (3) are similar.  If $\beta$ is $p$-restricted, then $\gamma^c = \beta'$ is $p$-regular, and this is equivalent to $\wt{\gamma}^c$ being $p$-regular, or $\wt{\beta} = (\wt{\gamma}^c)'$ being $p$-restricted.  Furthermore, we have $[S^{\delta^c} : D^{\gamma^c}] = 1 = [S^{\wt{\delta}^c} : D^{\wt{\gamma}^c}]$, so that part (2) follows by tensoring with the sign representation.
\end{proof}

Another block of interest in the study of the $[3:2]$-pair $(B, \wt{B})$ is the `intermediate' weight 4 block $\wb{B}$ of $\mathbb{F}\sym_{n-1}$.  The $p$-core of $\wb{B}$ can be obtained from that of $B$ by moving one bead from runner $i$ to runner $i-1$.  This block is `intermediate' in the following sense: if $M$ is a $B$-module, then $M \down{\wt{B}} = (M\down{\wb{B}})\down{\wt{B}}$, and if $N$ is a $\wt{B}$-module, then $N \up{B} = (N\up{\wb{B}})\up{B}$.

\begin{thm}[{\cite[Proposition 4.3]{mt2}}] \label{T:P}
Let $D^\lambda$ be a non-exceptional simple modules of $B$, and let $\wt{\lambda} = \Phi_{B,\wt{B}}(\lambda)$ and $\wb{\lambda} = \Phi_{B,\wb{B}}(\lambda)$.  Then
\begin{enumerate}
\item $D^{\wb{\lambda}} \up{B} = D^\lambda$, $D^{\wb{\lambda}} \down{\wt{B}} = D^{\wt{\lambda}}$;
\item $D^{\wt{\lambda}} \up{\wb{B}} \cong D^\lambda \down{\wb{B}}$, and $D^{\wt{\lambda}} \up{\wb{B}}$ is non-simple, has a simple head and a simple socle both isomorphic to $D^{\wb{\lambda}}$ and the composition factors $D^{\wb{\mu}}$ of its heart satisfy $D^{\wb{\mu}} \up{B} = 0 = D^{\wb{\mu}} \down{\wt{B}}$.
\end{enumerate}
\end{thm}

Let $\wb{\alpha} = \Phi_{B,\wb{B}}(\alpha)$.  Then we have the following multiplicities:
\begin{lem}[{\cite[Lemmas 4.5 and 4.7]{mt2}}] \label{L:ll3}
$[D^{\alpha} \down{\wb{B}} : D^{\wb{\alpha}}] = 3 = [D^{\wt{\alpha}} \up{\wb{B}} : D^{\wb{\alpha}}]$, and $[D^{\wb{\alpha}} \up{B} : D^{\alpha}] = 2 = [D^{\wb{\alpha}} \down{\wt{B}} : D^{\wt{\alpha}}]$.  In fact, $D^{\wb{\alpha}} \up{B}$ and $D^{\wb{\alpha}} \down{\wt{B}}$ have radical length $3$.
\end{lem}

\begin{rem}
It is worth mentioning that $D^{\wt{\alpha}} \up{\wb{B}} \not\cong D^\alpha \down{\wb{B}}$ (compare this with Theorem \ref{T:P}(2)).  This follows from the fact that the dimension of $\Hom_{\wb{B}}(D^{\alpha} \down{\wb{B}}, D^{\wt{\alpha}} \up{\wb{B}}) \cong \Hom_{\wt{B}}(D^{\alpha} \down{\wt{B}}, D^{\wt{\alpha}})$ is two, while $\End_{\wb{B}}(D^{\wt{\alpha}} \up{\wb{B}})$ is isomorphic to $\mathbb{F}[x]/(x^3)$ by Theorem \ref{T:Klesh} and Lemma \ref{L:ll3}, and hence has dimension three.
\end{rem}

From now on, we write $\Phi$ for $\Phi_{B, \wt{B}}$.  Also, following \cite{mt2}, write $L_5$ (resp.\ $\wt{L}_5$) for the indecomposable direct summand of $D^{\wt{\alpha}} \up{B}$ (resp.\ $D^\alpha \down{\wt{B}}$).  Thus, $D^{\wt{\alpha}} \up{B} \cong L_5^{\oplus 2}$ and $D^{\alpha} \down{\wt{B}} \cong \wt{L}_5^{\oplus 2}$.

\begin{prop} \label{P:phi}
There exists $\phi \in \End_B(P(D^{\wt{\alpha}}))$ such that $\phi(P(D^{\wt{\alpha}})) = \wt{L}_5$, $\phi(\wt{L}_5) = D^{\wb{\alpha}} \down{\wt{B}}$, $\phi(D^{\wb{\alpha}} \down{\wt{B}}) = D^{\wt{\alpha}}$.
\end{prop}

\begin{proof}
Since $P(D^{\wt{\alpha}}) \cong D^{\check{\alpha}} \up{\wt{B}}$ and $[P(D^{\wt{\alpha}}):D^{\wt{\alpha}}] =4$, we see that $\End_{\wt{B}}(P(D^{\wt{\alpha}})) \cong \mathbb{F}[x]/(x^4)$ by Theorem \ref{T:Klesh}.  Let $\phi \in \End_{\wt{B}}(P(D^{\wt{\alpha}}))$ such that $\{ 1, \phi, \phi^2, \phi^3 \}$ is a basis for $\End_{\wt{B}}(P(D^{\wt{\alpha}}))$ .  Then $P(D^{\wt{\alpha}})$, $\phi(P(D^{\wt{\alpha}}))$, $\phi^2(P(D^{\wt{\alpha}}))$ and $\phi^3(P(D^{\wt{\alpha}}))$ are submodules of $P(D^{\wt{\alpha}})$ with simple head $D^{\wt{\alpha}}$, with decreasing multiplicity of $D^{\wt{\alpha}}$ as a composition factor.  The Proposition thus follows from (the analogue of) \cite[Proposition 5.2]{mt2}.
\end{proof}

\begin{lem}[{\cite[Proposition 4.6]{mt2}}] \label{L:list}
Let $D^{\lambda}$ be a non-exceptional simple module of $B$, and let $\wt{\lambda} = \Phi(\lambda)$.  The following table provides all the possible composition multiplicities of $D^{\lambda}$ (resp.\ $D^{\wt{\lambda}}$) in $P(D^{\alpha})$, $L_5$ and $D^{\wb{\alpha}} \up{B}$ (resp.\ $P(D^{\wt{\alpha}})$, $\wt{L}_5$ and $D^{\wb{\alpha}} \down{\wt{B}}$):
$$
\begin{matrix}
& \text{\tiny{$[P(D^{\wt{\alpha}}) : D^{\wt{\lambda}}]$}} & \text{\tiny{$[\wt{L}_5 : D^{\wt{\lambda}}]$}} & \text{\tiny{$[D^{\wb{\alpha}} \down{\wt{B}} : D^{\wt{\lambda}}]$}} & \text{\tiny{$[D^{\wb{\alpha}} \up{B} : D^{\lambda}]$}} & \text{\tiny{$[L_5 : D^{\lambda}]$}} & \text{\tiny{$[P(D^{\alpha}) : D^{\lambda}]$}} \\[3pt]
& 0 & 0 & 0 & 0 & 0 & 0 \\
\textup{I} & 1 & 0 & 0 & 1 & 2 & 3 \\
\textup{II} & 2 & 1 & 0 & 0 & 1 & 2 \\
\textup{III} & 3 & 2 & 1 & 0 & 0 & 1 \\
\textup{IV(A)} & 4 & 2 & 1 & 1 & 2 & 4 \\
\textup{IV(B)} & 4 & 2 & 0 & 0 & 2 & 4
\end{matrix}
$$
\end{lem}

As a corollary, we have the following:

\begin{cor} \label{C:diffpty}
In cases \textup{I}, \textup{III} and \textup{IV(A)} of Lemma \ref{L:list}, we have $\pty \alpha \ne \pty \lambda$.
\end{cor}

\begin{proof}
In these cases, we have either $D^{\lambda}$ occurring as a composition factor of the heart of $D^{\wb{\alpha}} \up{B}$, or $D^{\Phi(\lambda)}$ occurring as a composition factor of the heart of $D^{\wb{\alpha}} \down{\wt{B}}$.  By Lemma \ref{L:ll3}, the hearts of $D^{\wb{\alpha}} \up{B}$ and $D^{\wb{\alpha}} \down{\wt{B}}$ are semi-simple.  Thus, $\pty{\wt{\lambda}} = \pty\lambda \ne \pty \alpha = \pty{\wt{\alpha}}$ in both instances by Theorem \ref{T:basicprop}(2) and Lemma \ref{L:pty}.
\end{proof}

\begin{prop} \label{P:samepty}
In cases \textup{II} and \textup{IV(B)} of Lemma \ref{L:list}, we have $\pty \alpha = \pty \lambda$.
\end{prop}

To prove Proposition \ref{P:samepty}, we use the $v$-decomposition numbers $d_{\lambda\mu}(v)$ arising from the canonical basis of the Fock space representation $\mathcal{F}$ of $U_v(\widehat{\mathfrak{sl}}_p)$.

Suppose that the fixed abacus display of the $p$-core of $B$ has $t$ beads.  Let $r$ be the residue class of $(i-t)$ modulo $p$.  In this section, we write $e$ and $f$ for the elements $e_r$ and $f_r$ of $U_v(\widehat{\mathfrak{sl}}_p)$ respectively.  Furthermore, $e^{(2)} = e^2/(v+v^{-1})$ and $f^{(2)} = f^2/(v+v^{-1})$.

\begin{lem} \label{L:e}
The table below records the effects of $e^{(2)}$ and $f^{(2)}$ on the standard basis elements of $\mathcal{F}$ labelled by exceptional partitions in $B$ and $\wt{B}$ respectively.  The entry on a row labelled by $\rho$ and a column labelled by $\wt{\sigma}$ is $\langle e^{(2)} s(\rho), s(\wt{\sigma}) \rangle$ ($=\langle s(\rho), f^{(2)} s(\wt{\sigma}) \rangle$).
$$
\begin{matrix}
\quad\ & \wt{\alpha} & \wt{\beta} & \wt{\gamma} & \wt{\delta}\\[6pt]
\alpha & v^{-2} & v^{-1} & 1 & 0 \\[3pt]
\beta  & v^{-1} & 1      & 0 & 1 \\[3pt]
\gamma & 1      & 0      & 1 & v \\[3pt]
\delta & 0      & 1      & v & v^2
\end{matrix}
$$
\end{lem}

\begin{lem} \label{L:e2}
$G(\alpha) = s(\alpha) + vs(\beta) + v^2 s(\gamma) + v^3s(\delta)$ and $G(\wt{\alpha}) = s(\wt{\alpha}) + vs(\wt{\beta}) + v^2 s(\wt{\gamma}) + v^3s(\wt{\delta})$.
\end{lem}

\begin{proof}
Since $\check{\alpha}$ has weight $0$, we see that $G(\check{\alpha}) = s(\check{\alpha})$ (\cite[Theorem 6.8(iii)]{LLT}).  This gives $fG(\check{\alpha}) = fs(\check{\alpha}) = s(\wt{\alpha}) + vs(\wt{\beta}) + v^2 s(\wt{\gamma}) + v^3s(\wt{\delta})$.  But this implies that $fG(\check{\alpha}) = G(\wt{\alpha})$ (\cite[Theorem 6.1]{LLT}) so that the second assertion follows.  Similar arguments apply to the first assertion.
\end{proof}

\begin{prop} \label{P:caseII}
Suppose $\lambda$ is in case \textup{II} of Lemma \ref{L:list}, and let $\wt{\lambda} = \Phi(\lambda)$.
The following is a complete list of possible $v$-decomposition numbers $d_{\mu\lambda}(v)$, $d_{\wt{\mu}\wt{\lambda}}(v)$ where $\mu$ are $\wt{\mu}$ are exceptional partitions in $B$ and $\wt{B}$ respectively.
$$
\begin{matrix}
d_{\alpha\lambda}(v) & d_{\beta\lambda}(v) & d_{\gamma\lambda}(v) & d_{\delta\lambda}(v) &
d_{\wt{\alpha}\wt{\lambda}}(v) & d_{\wt{\beta}\wt{\lambda}}(v) & d_{\wt{\gamma}\wt{\lambda}}(v) & d_{\wt{\delta}\wt{\lambda}}(v) \\
v^2 & v^3 & \cdot & \cdot & v^2 & v^3 & \cdot & \cdot \\
v^2 & \cdot & v^2 & \cdot & v^2 & \cdot & v^2 & \cdot \\
v^2 & \cdot & \cdot & v & \cdot & v & v^2 & \cdot \\
\cdot & v & v^2 & \cdot & v^2 & \cdot & \cdot & v \\
\cdot & v & \cdot & v & \cdot & v & \cdot & v \\
\cdot & \cdot & 1 & v & \cdot & \cdot & 1 & v
\end{matrix}
$$
\end{prop}

\begin{proof}
By Lemma \ref{L:list}, we have $[\wt{L}_5: D^{\wt{\lambda}}] = 1$.  Thus, $P(D^{\wt{\lambda}}) \up{B} \cong P(D^{\lambda})^{\oplus 2} \oplus P(D^{\alpha})^{\oplus 2}$ by Frobenius reciprocity.  These arguments carry over to the Iwahori-Hecke algebras in characteristic zero, so that $P(\mathcal{D}^{\wt{\lambda}}) \up{B} \cong P(\mathcal{D}^{\lambda})^{\oplus 2} \oplus P(\mathcal{D}^{\alpha})^{\oplus 2}$, and hence $f^{(2)} G(\wt{\lambda}) = G(\lambda) + G(\alpha)$ by Proposition \ref{P:Ariki} (and the Note immediately after that).  Thus,
$$
d_{\mu\lambda}(v) + d_{\mu\alpha}(v) = \langle G(\lambda) + G(\alpha), s(\mu) \rangle =
\langle f^{(2)} G(\wt{\lambda}), s(\mu) \rangle
= \langle G(\wt{\lambda}), e^{(2)} s(\mu) \rangle
$$
for all partitions $\mu$.  Varying $\mu$ over the exceptional partitions in $B$ and using Lemmas \ref{L:e} and \ref{L:e2}, we obtain the following four equations:
$$
\begin{array}{llllrrrrrrr}
d_{\alpha\lambda}(v)&+ &1 &= &v^{-2} d_{\wt{\alpha}\wt{\lambda}}(v) &+ &v^{-1} d_{\wt{\beta}\wt{\lambda}}(v) &+  &d_{\wt{\gamma}\wt{\lambda}}(v) & & \\[4pt]
d_{\beta\lambda}(v) &+ &v &= &v^{-1} d_{\wt{\alpha}\wt{\lambda}}(v) &+ &d_{\wt{\beta}\wt{\lambda}}(v) &+ & &  &d_{\wt{\delta}\wt{\lambda}}(v) \\[4pt]
d_{\gamma\lambda}(v) &+ &v^2 &= &d_{\wt{\alpha}\wt{\lambda}}(v) &+ &&&d_{\wt{\gamma}\wt{\lambda}}(v) &+ &v d_{\wt{\delta}\wt{\lambda}}(v) \\[4pt]
d_{\delta\lambda}(v) &+ &v^3 &= &&&d_{\wt{\beta}\wt{\lambda}}(v) &+ &v d_{\wt{\gamma}\wt{\lambda}}(v) &+  &v^2 d_{\wt{\delta}\wt{\lambda}}(v)
\end{array}
$$
As exactly two among $d_{\wt{\alpha}\wt{\lambda}}(v)$, $d_{\wt{\beta}\wt{\lambda}}(v)$, $d_{\wt{\gamma}\wt{\lambda}}(v)$ and $d_{\wt{\delta}\wt{\lambda}}(v)$ are zero, while the other two are monic monomials (since $[P(D^{\wt{\alpha}}):D^{\wt{\lambda}}] = 2$), the above equations give all the possibilities as listed in the Proposition.  We illustrate this with an example.  Suppose $d_{\wt{\gamma}\wt{\lambda}}(v) =d_{\wt{\delta}\wt{\lambda}}(v) = 0$.  Then the last two equations give $d_{\wt{\alpha}\wt{\lambda}}(v) = v^2$, $d_{\wt{\beta}\wt{\lambda}}(v) = v^3$, $d_{\delta\lambda}(v) = d_{\gamma\lambda}(v) = 0$.  Substituting these into the first two equations, we get $d_{\alpha\lambda}(v) = v^2$, $d_{\beta\lambda}(v) = v^3$.  This is precisely the first row in the Proposition.
\end{proof}

\begin{prop} \label{P:caseIV}
Suppose $[P(D^{\alpha}) : D^{\lambda}] =4$ with $\lambda \ne \alpha$, and let $\wt{\lambda} = \Phi(\lambda)$.
The following is a complete list of possible $v$-decomposition numbers $d_{\mu\lambda}(v)$, $d_{\wt{\mu}\wt{\lambda}}(v)$ where $\mu$ and $\wt{\mu}$ are exceptional partitions in $B$ and $\wt{B}$ respectively.
$$
\begin{matrix}
d_{\alpha\lambda}(v) & d_{\beta\lambda}(v) & d_{\gamma\lambda}(v) & d_{\delta\lambda}(v) &
d_{\wt{\alpha}\wt{\lambda}}(v) & d_{\wt{\beta}\wt{\lambda}}(v) & d_{\wt{\gamma}\wt{\lambda}}(v) & d_{\wt{\delta}\wt{\lambda}}(v) \\
v & v^2 & v & v^2 & v & v^2 & v & v^2 \\
v^2 & v & v^2 & v & v^2 & v & v^2 & v
\end{matrix}
$$
\end{prop}

\begin{proof}
Note that $\lambda \notin \{ \alpha, \beta, \gamma, \delta, m(\alpha)', m(\beta)', m(\gamma)', m(\delta)' \}$.  Thus if $\mu$ is an exceptional partition in $B$, then
$$
d_{\mu\lambda}(v) =
\begin{cases}
v, &\text{if } \pty \lambda \ne \pty\mu; \\
v^2, & \text{if } \pty \lambda = \pty\mu
\end{cases}
$$
by Corollary \ref{C:monomials}.  As analogous statements hold for $d_{\wt{\mu}\wt{\lambda}}(v)$, the Proposition follows immediately from Lemma \ref{L:pty}.
\end{proof}

\begin{proof}[Proof of Proposition \ref{P:samepty}]
By Proposition \ref{P:caseII} and Corollary \ref{C:monomials}, we see that $\pty\lambda = \pty\alpha$ if $\lambda$ is in case II of Lemma \ref{L:list}.  By Proposition \ref{P:caseIV} and Corollary \ref{C:monomials}, it suffices to show that $d_{\wt{\alpha}\wt{\lambda}}(v) \ne v$, where $\wt{\lambda} = \Phi(\lambda)$, when $\lambda$ is in case IV(B) of Lemma \ref{L:list}.  Indeed, in this case, we have $[D^{\wb{\alpha}} \down{\wt{B}}: D^{\wt{\lambda}}] = 0$ by Lemma \ref{L:list}, so that $P(D^{\wt{\lambda}}) \up{\wb{B}} \cong P(D^{\wb{\lambda}})$ by Frobenius reciprocity and Theorem \ref{T:P}(1), where $\wb{\lambda} = \Phi_{B,\wb{B}}(\lambda)$.  These arguments carry over to the Iwahori-Hecke algebra in characteristic zero, so that $P(\mathcal{D}^{\wt{\lambda}}) \up{\wb{B}} \cong P(\mathcal{D}^{\wb{\lambda}})$, and hence $fG(\wt{\lambda}) = G(\wb{\lambda})$ by Proposition \ref{P:Ariki} (and the Note immediately after that).  Thus $\langle fG(\wt{\lambda}), s(\wb{\alpha}) \rangle \in v\mathbb{N}_0[v]$.  This implies that $d_{\wt{\alpha}\wt{\lambda}}(v)\langle f s(\wt{\alpha}), s(\wb{\alpha}) \rangle = v^{-1} d_{\wt{\alpha}\wt{\lambda}}(v) \in v \mathbb{N}_0[v]$, so that $d_{\wt{\alpha}\wt{\lambda}}(v) \ne v$.
\end{proof}

\begin{rem}
The sufficient condition (Y3) in \cite{mt2} holds trivially by Proposition \ref{P:samepty} and Theorem \ref{T:basicprop}(2).
\end{rem}

With Proposition \ref{P:samepty}, we conclude that $\Ext^1(D^{\lambda},D^{\alpha}) = 0 = \Ext^1(D^{\wt{\lambda}},D^{\wt{\alpha}})$ when $\lambda$ is in Cases II or IV(B).  We are also able to conclude that $\Ext^1(D^{\wt{\lambda}},D^{\wt{\alpha}})$ is non-zero when $\lambda$ is in Cases III or IV(A) since in these cases, $D^{\wt{\lambda}}$ lies in the semi-simple heart of $D^{\wb{\alpha}} \down{\wt{B}}$.  However, it is as yet unclear if $\Ext^1(D^{\wt{\lambda}},D^{\wt{\alpha}})$ is non-zero when $\lambda$ is in Case I, even though we know that $\pty\lambda \ne \pty\alpha$ by Corollary \ref{C:diffpty}.  We now address this issue.

\begin{lem} \label{L:1}
Suppose $[P(D^{\wt{\alpha}}):D^{\wt{\lambda}}] = 1 = [S^{\wt{\alpha}}:D^{\wt{\lambda}}]$.  Then $\wt{\alpha}$ is $p$-restricted and $\wt{\alpha} = m(\wt{\lambda})'$.
\end{lem}

\begin{proof}
By Lemma \ref{L:list}, we have $[\wt{L}_5: D^{\wt{\lambda}}] = 0$.  Thus, $P(D^{\wt{\lambda}}) \up{B} \cong P(D^{\lambda})^{\oplus 2}$ by Frobenius reciprocity.   These arguments carry over to the Iwahori-Hecke algebra in characteristic zero, so that $P(\mathcal{D}^{\wt{\lambda}}) \up{B} \cong P(\mathcal{D}^{\lambda})^{\oplus 2}$, and hence $f^{(2)} G(\wt{\lambda}) = G(\lambda)$ by Proposition \ref{P:Ariki} (and the Note immediately after that).  Thus,
$$
v^{-2}d_{\wt{\alpha}\wt{\lambda}}(v) = \langle f^{(2)} G(\wt{\lambda}), s(\alpha) \rangle = \langle G(\lambda), s(\alpha) \rangle \in v\mathbb{N}_0[v].$$
This implies that $d_{\wt{\alpha}\wt{\lambda}}(v) = v^3$ and $\wt{\alpha} = m(\wt{\lambda})'$ by Corollary \ref{C:monomials}.
\end{proof}

\begin{prop} \label{P:1}
Suppose $[P(D^{\wt{\alpha}}):D^{\wt{\lambda}}] = 1$.  Then $\Ext^1(D^{\wt{\lambda}},D^{\wt{\alpha}}) = 0$.
\end{prop}

\begin{proof}
Suppose for the sake of contradiction that $\Ext^1(D^{\wt{\lambda}},D^{\wt{\alpha}}) \ne 0$.
If $\wt{\beta}$ is $p$-regular, then $D^{\wt{\lambda}}$ lies in the second radical layer of $S^{\wt{\alpha}}$ by Corollary 5.6 of \cite{mt2}.  However, by Lemma \ref{L:1}, we also have $\wt{\alpha}$ is $p$-restricted and $\wt{\alpha} = m(\wt{\lambda})'$, so that $D^{\wt{\lambda}}$ occurs as the socle of $S^{\wt{\alpha}}$.  This is impossible, since $S^{\wt{\alpha}}$ has a nonzero heart, with composition factor such as $m(\wt{\beta}')$ by Lemma \ref{L:conjugate} (note that $\wt{\beta}$ is $p$-restricted since $\wt{\alpha}$ is).

If $\wt{\gamma}$ is $p$-restricted, then applying the argument in the last paragraph to the conjugate blocks of $(B,\wt{B})$, we obtain $\Ext^1(D^{m(\wt{\lambda})},D^{m(\wt{\alpha})}) \ne 0$, so that $\Ext^1(D^{\wt{\lambda}},D^{\wt{\alpha}}) \ne 0$.

Since we cannot have $\wt{\beta}$ being $p$-singular and $\wt{\gamma}$ being non-$p$-restricted at the same time, we are done.
\end{proof}

\begin{rem} \hfill
\begin{enumerate}
\item The analogue of Proposition \ref{P:1} also holds --- if $[P(D^{\alpha}):D^{\lambda}] = 1$, then $\Ext^1(D^{\lambda},D^{\alpha}) = 0$.  Its proof however is more complicated, as $\alpha$ being $p$-restricted no longer implies $\beta$ (or $\gamma$) being $p$-restricted; also it is possible for $\beta$ to be $p$-singular and $\gamma$ to be non-$p$-restricted at the same time.
\item By Proposition \ref{P:1} and its analogue, it is straightforward to see that the sufficient condition (Y2) in \cite{mt2} holds.
\end{enumerate}
\end{rem}

\begin{prop} \label{P:L_5}
The radical length of $\wt{L}_5$ is $5$.
\end{prop}

\begin{proof}
Let $\wt{M}$ be the submodule of $\wt{L}_5$ such that $\wt{L}_5/\wt{M} \cong D^{\wb{\alpha}} \down{\wt{B}}$.  Since $\wt{M}$ is the largest submodule of $\wt{L}_5$ with $[\wt{M} : D^{\wt{\alpha}}] = 1$, and $[\rad(D^{\wb{\alpha}} \down{\wt{B}}) : D^{\wt{\alpha}}] = 1$, we have $\rad(D^{\wb{\alpha}} \down{\wt{B}}) \subseteq \wt{M}$.
By Lemma \ref{L:list}, the composition factors of $\wt{M}/\rad(D^{\wb{\alpha}} \down{\wt{B}})$ are in cases II or IV(B), and hence are indexed by partitions having the same relative sign as $\wt{\alpha}$ by Proposition \ref{P:samepty} and Lemma \ref{L:pty}.  Thus, these composition factors lie in an odd radical layer of $\wt{L}_5$, and hence in the third radical layer, since $\wt{L}_5/\wt{M} \cong D^{\wb{\alpha}} \down{\wt{B}}$ has radical length $3$ by Lemma \ref{L:ll3}.  This implies that $\wt{L}_5$ has radical length $5$.
\end{proof}

\begin{thm} \label{T:alpha}
The radical length of $P(D^{\wt{\alpha}})$ is $7$.
\end{thm}

\begin{proof}
Let $\phi \in \End_{\wt{B}}(P(D^{\wt{\alpha}}))$ be as described in Proposition \ref{P:phi}.  Then
$$
\rad(D^{\wb{\alpha}} \down{\wt{B}}) = \ker(\phi) \cap D^{\wb{\alpha}} \down{\wt{B}}
= \ker(\phi|_{\wt{L}_5}) \cap D^{\wb{\alpha}} \down{\wt{B}}.
$$
Let $M = D^{\wb{\alpha}} \down{\wt{B}} + \ker(\phi)$, and $N = D^{\wb{\alpha}} \down{\wt{B}} + (\ker(\phi) \cap \wt{L}_5) = D^{\wb{\alpha}} \down{\wt{B}} + \ker(\phi|_{\wt{L}_5})$,  and consider the filtration
$$
0 \subseteq D^{\wb{\alpha}} \down{\wt{B}} \subseteq N \subseteq M \subseteq P(D^{\wt{\alpha}}).
$$
We have $P(D^{\wt{\alpha}})/M \cong (P(D^{\wt{\alpha}})/\ker(\phi)) / (M/\ker(\phi))$, and $P(D^{\wt{\alpha}})/\ker(\phi) \cong \phi(P(D^{\wt{\alpha}})) = \wt{L}_5$ while
$$M/\ker(\phi) \cong D^{\wb{\alpha}} \down{\wt{B}}/(D^{\wb{\alpha}} \down{\wt{B}} \cap \ker(\phi)) =  D^{\wb{\alpha}} \down{\wt{B}} / \rad(D^{\wb{\alpha}} \down{\wt{B}}) = D^{\wt{\alpha}},$$
so that $P(D^{\wt{\alpha}})/M \cong \wt{L}_5/ D^{\wt{\alpha}}$.  Thus $P(D^{\wt{\alpha}})/M$ has radical length $4$ by Proposition \ref{P:L_5}.  Next, $M/N = (N + \ker(\phi))/N \cong \ker(\phi)/(N \cap \ker(\phi)) = \ker(\phi)/ (\ker(\phi|_{\wt{L}_5}) + (D^{\wb{\alpha}} \up{B} \cap \ker(\phi))) = \ker(\phi)/\ker(\phi|_{\wt{L}_5})$.
By Lemma \ref{L:list}, we see that the composition factors of $M/N$ are in cases I or IV(A), and are thus labelled by partitions having the same relative sign, which is different from that of $\wt{\alpha}$ by Corollary \ref{C:diffpty} and Lemma \ref{L:pty}.  Thus they lie in even radical layers of $P(D^{\wt{\alpha}})$ by Theorem \ref{T:basicprop}(2).  Since $P(D^{\wt{\alpha}})/M$ has radical length $4$, this implies that $P(D^{\wt{\alpha}})/N$ also has radical length $4$.  Now, as $N \subseteq \rad(\wt{L}_5)$, and hence has radical length at most $4$ by Proposition \ref{P:L_5}, we see that $P(D^{\wt{\alpha}})$ has radical length at most $8$.  But at the same time, the radical length of $P(D^{\wt{\alpha}})$ is greater than $5$ (since $L_5$ is a proper submodule and has radical length $5$ by Proposition \ref{P:L_5}), and must be odd (since $D^{\alpha}$ and hence $P(D^{\alpha})$ are self-dual and the $\Ext$-quiver of $\wt{B}$ is bipartite).  The Theorem thus follows.
\end{proof}

\begin{rem}
In view of the symmetry between $B$ and $\wt{B}$, statements and proofs entirely analogous to Proposition \ref{P:L_5} and Theorem \ref{T:alpha} also hold for the block $B$.
\end{rem}

\begin{cor}[of proof] \label{C:radlengkerphi}
Let $\phi$ be as described in Proposition \ref{P:phi}.  The radical length of $\ker(\phi)$ is at most $4$.
\end{cor}

\begin{proof}
Keeping the notations used in the proof of Theorem \ref{T:alpha}, we have proved the following:
\begin{itemize}
\item The radical length of $N$ is at most $4$;
\item The composition factors of $M/N$ are labelled by partitions all having the same relative sign, which is different from that of $\wt{\alpha}$.
\end{itemize}
This shows that the radical length of $M$ is at most $4$ too, since $M$ (as well as $N$) has a simple socle $D^{\wt{\alpha}}$.  Thus, the same holds for $\ker(\phi)$ since $\ker(\phi)$ is a submodule of $M$.
\end{proof}

\begin{prop} \label{P:geq7}
Let $D^{\wt{\lambda}}$ be a non-exceptional simple module of $\wt{B}$ and suppose $[P(D^{\wt{\alpha}}): D^{\wt{\lambda}}] > 0$.  Then the radical length of $P(D^{\wt{\lambda}})$ is at least $7$.
\end{prop}

\begin{proof}
We show $\rad^6 (P(D^{\wt{\lambda}})) \ne 0$.  Note that if the radical length of $P(D^{\wt{\lambda}})$ is $\ell$, then $\ell$ is odd by Theorem \ref{T:basicprop}(2), and $\rad^{\ell-1}(P(D^{\wt{\lambda}})) = D^{\wt{\lambda}}$.  We consider the various cases of Lemma \ref{L:list} separately.
\begin{description}
\item[Cases III and IV(A)]
For these cases, $D^{\wt{\lambda}}$ is a composition factor of $\rad(D^{\wb{\alpha}} \down{\wt{B}})$ by Lemma \ref{L:list}.  Since $D^{\wb{\alpha}} \down{\wt{B}} \subseteq \rad^2(\wt{L}_5)$ and $\wt{L}_5 \subseteq \rad^2(P(D^{\wt{\alpha}}))$, we see that $D^{\wt{\lambda}}$ is a composition factor of $\rad^5(P(D^{\wt{\alpha}}))$.  As the simple modules of symmetric groups are self-dual, this implies that $D^{\wt{\alpha}}$ is a composition factor of $\rad^5(P(D^{\wt{\lambda}}))$.  Thus $\rad^6(P(D^{\wt{\lambda}})) \ne 0$.
\item[Cases II and IV(B)]
For these cases, $D^{\wt{\lambda}}$ is a composition factor of $\wt{L}_5$ by Lemma \ref{L:list}.  Since $\pty{\wt{\lambda}} = \pty{\wt{\alpha}}$ for these cases by Proposition \ref{P:samepty}, we see that in fact $D^{\wt{\lambda}}$ is a composition factor of $\rad^2(\wt{L}_5)$, and hence of $\rad^4(P(D^{\wt{\alpha}}))$.  By self-duality of simple modules, this implies that $D^{\wt{\alpha}}$ is a composition factor of $\rad^4(P(D^{\wt{\lambda}}))$.  Thus $\rad^6(P(D^{\wt{\lambda}})) \ne 0$.
\item[Case I]  Since $\pty{\wt{\alpha}} \ne \pty{\wt{\lambda}}$ by Corollary \ref{C:diffpty}, and $\Ext^1(D^{\wt{\lambda}},D^{\wt{\alpha}}) = 0$ by Proposition \ref{P:1}, $D^{\wt{\lambda}}$ is a composition factor of $\rad^3(P(D^{\wt{\lambda}}))$.  By self-duality of simple modules, this implies that $D^{\wt{\alpha}}$ is a composition factor of $\rad^3(P(D^{\wt{\lambda}}))$.  Now if $\rad^6(P(D^{\wt{\lambda}})) = 0$, then $\rad^4(P(D^{\wt{\lambda}})) = D^{\wt{\lambda}}$, so that $D^{\wt{\alpha}}$ extends $D^{\wt{\lambda}}$, a contradiction.
\end{description}
\end{proof}

The next result is on general representation theory of finite groups and will be needed in the proof of Proposition \ref{P:32}.

\begin{prop} \label{P:pre}
Let $H$ be a subgroup of a finite group $G$, and let $B$ be a block of $\mathbb{F}G$ and $C$ be a block of $\mathbb{F}H$.  Suppose $N$ is a submodule of a $C$-module $M$ such that for any simple $B$-module $S$, and any $\phi \in \Hom_{C}(M, S \down{C})$, we have $N \subseteq \ker(\phi)$.  Then $N \up{B} \subseteq \rad(M \up{B})$.
\end{prop}

\begin{proof}
This follows from Frobenius reciprocity.  Suppose $N \up{B} \nsubseteq \rad(M \up{B})$.  Then there exists a simple $B$-module $S$ and a map $\psi \in \Hom_{\mathbb{F}G}(M \up{G}, S)$ such that $\psi(N \up{G}) \ne 0$.  Thus $\psi(1 \otimes n) \ne 0$ for some $1 \otimes n \in \mathbb{F}G \otimes_{\mathbb{F}H} N = N \up{G}$.  Define $\phi : M \to S\down{C}$ by $\phi(m) = e_C\psi(1 \otimes m)$, where $e_C$ is the block idempotent of $C$.  Then $\phi \in \Hom_{C}(M, S\down{C})$, and
$$
\phi(n) = e_C\psi(1 \otimes n) = \psi(1 \otimes e_Cn) = \psi(1 \otimes n) \ne 0,
$$
a contradiction.
\end{proof}

We conclude this paper with a proof of Proposition \ref{P:32}.

\begin{proof}[Proof of Proposition \ref{P:32}]
Let $D^{\wt{\lambda}}$ be a simple module of $\wt{B}$, and let $l$ be the radical length of $P(D^{\wt{\lambda}})$.  If $[P(D^{\wt{\lambda}}) : D^{\wt{\alpha}} ] =0$, then $P(D^{\wt{\lambda}})$ has the same radical length as $P(D^{\wt{\lambda}}) \up{B}$, so that $l = 7$.  We thus assume that $[P(D^{\wt{\lambda}}) : D^{\wt{\alpha}} ] > 0$, and we consider the cases of $\wt{\lambda}$ having the same and different relative sign as $\wt{\alpha}$ separately.
\begin{description}
\item[Case 1. $\pty{\wt{\lambda}} = \pty{\wt{\alpha}}$] \hfill
\begin{description}
\item[Subcase 1a. $\wt{\lambda} = \wt{\alpha}$]  This is dealt with in Theorem \ref{T:alpha}.
\item[Subcase 1b. $\wt{\lambda} \ne \wt{\alpha}$]  By Proposition \ref{P:geq7}, $l \geq 7$.  Suppose that $l > 7$; so we have $l \geq 9$.  Let $M$ be a submodule of $\rad^2(P(D^{\wt{\lambda}}))$ having a simple head, say $D^{\wt{\mu}}$, and radical length $l-2$.  Note that $\wt{\mu} \ne \wt{\alpha}$, since $P(D^{\wt{\alpha}})$ has radical length $7$.  Let $N$ be a submodule of $\rad^2(M)$ having a simple head, say $D^{\wt{\nu}}$, and radical length $l-4$.  Clearly, $[N : D^{\wt{\alpha}}] \leq 1$ with equality if and only if $\wt{\nu} = \wt{\alpha}$.  Suppose for a contradiction that $\wt{\nu} = \wt{\alpha}$.  Embed $N^*$ (the dual of $N$) into $P(D^{\wt{\alpha}})$.  Let $\phi \in \End_{\wt{B}}(P(D^{\wt{\alpha}}))$ be as described in Proposition \ref{P:phi}.   Then $N^* \subseteq \ker(\phi)$ since $\ker(\phi)$ is the unique maximal submodule of $P(D^{\alpha})$ such that $[\ker(\phi) : D^{\alpha}] \leq 1$.  But $\ker(\phi)$ has radical length at most $4$ by Corollary \ref{C:radlengkerphi}, contradicting $N^*$ having radical length $l-4$ ($\geq 5$).  Thus, $\wt{\nu} \ne \wt{\alpha}$ and hence $[N : D^{\wt{\alpha}}] =0$.  This implies that $N \up{B}$ has radical length $l-4$.  Note that $M$ and $N$ satisfy the hypothesis of Proposition \ref{P:pre}; thus $N \up{B} \subseteq \rad(M \up{B})$.  Furthermore, by Frobenius reciprocity, only $D^{\mu}$, where $\mu = \Phi^{-1}_{B,\wt{B}} (\wt{\mu})$, and possibly $D^{\alpha}$ occur in the head of $M \up {B}$, while the head of $N \up{B}$ is isomorphic to $D^{\wt{\nu}} \up{B}$ ($\cong (D^{\nu})^{\oplus 2}$ say).  Since $\pty \mu = \pty \nu = \pty\alpha$, we see that in fact $N \up{B} \subseteq \rad^2(M \up{B})$, so that $M \up{B}$ has radical length at least $l-2$.  Since $P(D^{\wt{\lambda}})$ and $M$ also satisfy the hypothesis of Proposition \ref{P:pre}, we have $M \up{B} \subseteq \rad(P(D^{\wt{\lambda}}) \up{B})$, so that $M \up{B}$ has radical length at most $6$.  This gives $l-2 \leq 6$, a contradiction.
\end{description}
\item[Case 2. $\pty{\wt{\lambda}} \ne \pty{\wt{\alpha}}$]  Let $P$ be the projective cover of $\rad(P(D^{\wt{\lambda}}))$.  Then $P$ is a direct sum of projective indecomposable modules which are all indexed by partitions having the same relative sign as $\wt{\alpha}$.  From Case (1), we conclude that $P$ has radical length $7$, so that $\rad(P(D^{\wt{\lambda}}))$ has radical length at most $6$.  Thus $l \leq 7$, and hence $l=7$ by Proposition \ref{P:geq7}.
\end{description}
\end{proof}

\end{document}